\documentclass{ctexart}
\usepackage{amsmath}
\usepackage{amsfonts}

\usepackage{geometry}                % See geometry.pdf to learn the layout options. There are lots.
\usepackage[parfill]{parskip}    % Activate to begin paragraphs with an empty line ra than an indent
\usepackage{graphicx}
\usepackage{amssymb}
\usepackage{epstopdf}
\usepackage[all]{xy}

\usepackage{amsthm} 
\theoremstyle{definition} 
\newtheorem{dfn}{Definition}[section] 
\theoremstyle{axiom}

\theoremstyle{remark} 
 
\theoremstyle{plain}

\theoremstyle{plain}

\usepackage{amssymb, wasysym, animate, multicol}

\usepackage{fontspec}
\usepackage{polyglossia}
\usepackage{xcolor}

\setmainlanguage{english}
\setotherlanguages{sanskrit} %% or other languages
\begin{document}

\title{From Schwartz Space to Structural Aging\\
(Towards a small history of Old Age mathematics)}

\author{Daniel Parrochia}
\date{University of Lyon (France)}
\maketitle

\vspace{1\baselineskip}
\hspace{8.7\baselineskip}
\textit{To the young, who will grow old, and to those who already are.} \\

\textbf{Abstract.}
Starting with some historical considerations, we examine the various attempts to model old age that fall within the realm of mathematical physics. For example, we study the mathematical theory of declining functions, which, since Laurent Schwartz's work on distributions, are generally defined on what is called a "Schwartz space." We explain how this essentially analytic formalization applies to the representation of the decline of vitality (or of a number of vital functions). We then show that there exist many other possible formalizations that utilize the resources of algebra, in particular the theory of "inclines" of Cao, Kim, and Roush and that of relational structures, on which an algebraic notion of "age," in Cameron's sense, can be defined. But we maintain, in conclusion, that a more positive representation of old age, consistent, moreover, with what a number of cultural sources from Antiquity to the 20th century tell us, should mobilize another type of formalization, capable of selecting, within the overall debate between the living being and its environment, persistent structures that do not decline amidst the decline of others.\\


\section{A Brief History of Old Age}
An informed study of Greek and Roman conceptions of old age reveals that, for ancient philosophers, the quality and meaning of this stage of life derive primarily from the individual's wisdom and moral character. The acknowledged disagreements among these philosophers often stem from divergences in their metaphysical, epistemological, and -- naturally -- moral doctrines. For instance, Plato suggests that old age is synonymous with liberation from the earthly desires of youth -- desires that merely distract man from truth, knowledge, and the good life (see {\it Rep.} IV, 428e–429a, VII, 540a–c, VIII, 563a7–b2; and {\it Laws} III, 690a–b and 714e)\footnote{Although political leaders (guardians, lawmakers) are generally chosen from among the middle-aged, Plato does not appear to posit a systematic link between old age and wisdom. Indeed, mindful of the infirmities of old age, he notes -- contrary to a famous adage by Solon -- that it is easier for an old man to run than to learn something new ({\it Rep.} VII, 536d).}. While Aristotle’s view is, on this subject, more negative\footnote{For him, old age is associated with physical decline ({\it phtisis}) ({\it De Mem.} 2, 453b5-6; {\it De Resp.} 10, 475b16-19; 476a7-8; {\it De Caelo}, II, 7, 288b13-17) -- most often accompanied by moral and intellectual degeneration, as illustrated by the portrait in {\it Rhet.} II, 13.}, his moral psychology acknowledges that the virtuous elderly person has many reasons to rejoice. Since virtue is self-perpetuating, one who acquires it before old age should be able to endure hardships with grace and fully savor contemplation in later life. Cicero, Seneca, and Plutarch\footnote{We shall return later to Cicero’s {\it De Senectute} (see section 11). Regarding Plutarch, his treatise {\it On Moral Virtue} (450 F1-7) may have influenced Plato’s conceptions. For Stoic views, see Seneca, {\it Letters to Lucilius}, 12 and 26; and also: {\it On the Shortness of Life}.} all maintain that the activities characteristic of old age are superior to any other way of life (see \cite{Mon}; \cite{Ant}). As is well known, the Middle Ages largely relied on the legacy of Antiquity; for a long time, old age was rarely treated as a subject in its own right, but rather viewed within the broader context of the "ages of life." These were typically categorized -- following the Pythagorean model -- into three or four periods (childhood, adolescence, maturity, old age), as seen in the works of Galen or Philippe de Novare, or sometimes more (six periods according to Isidore of Seville, seven according to Augustine and Montaigne); in any case, the final stage was invariably a disaster.

For everything gradually fades away, as Jaucourt explains in the entry on "Old Age" in Diderot and d'Alembert's {\it Encyclopédie}: physically, "as one advances in age, the bones, cartilages, membranes, flesh, and all the body's fibers grow dry and rigid"; fluids degenerate, arteries thicken, and "muscles lose their spring; the head wavers, the hand trembles, the legs falter; hearing, sight, and smell weaken, and even the sense of touch grows dull." The situation is obviously even worse for women, and the mental or moral aspect is no better: "One grows old in years without gaining an inch of wisdom," as Montaigne noted, and habit shackles the elderly just as surely as stiffness affects their movements. In the "Philosophy of Mind" section of the {\it Shorter Encyclopaedia of the Philosophical Sciences}, Hegel appears as one of the last heirs to this perspective: "An elderly person lives without any specific interest, for they have renounced the hope of realizing the ideals they once cherished, and the future seems to offer them neither promise nor novelty. On the contrary, they already regard themselves as familiar with the universal -- the essential principle underlying everything they might yet encounter" (see \cite{Pet}, vol. 2, § 396, Addition). In short, Geras (the god of Old Age), a son of Erebus and Night in mythology, is a terrible and implacable force for the elderly individual, compelling them to dissolve into a distant impersonality.

However, one must exercise particular caution regarding the subject of old age, for a preliminary question must be asked before any further reflection: what kind of old age are we talking about? In France, life expectancy at birth -- measured annually for the entire population since 1740 -- was still very low in the mid-18th century. Following successive wars, it stood at around 38.0 years between 1650 and 1679, rising to 42.9 years between 1770 and 1819\footnote{Molière's famous "old fogies" were thus, at most, around fifty years old.}. Although slightly higher for nobles than for the general population, it remained far lower for the latter than life expectancy at age 25 in 2020 (54.8 years for men, 60.6 years for women) (see \cite{Ins}). Based on mortality rates by age and sex from 1730–1749, 24\% of nobles and 60\% of the total population died before the age of 25. Yet since the early 19th century--a period when French life expectancy was at its lowest, averaging 33 years (see \cite{Clo}) -- it has risen steadily. Between 1900 and 2000, it climbed from 48 to 79 years, an increase of 65\% in a single century. Consequently, the questions raised by the process of aging stem essentially from a recent area of ​​inquiry\footnote{As late as 1983, Philippe Ariès lamented the lack of a true history of old age (see \cite{Ari}). That is, however, no longer the case. In fact, as early as 1982, some thought had already been given to the issue (see \cite{Roy}). Subsequently, however, the topic flourished (see, for example, \cite{Min}, \cite{Boi}, \cite{Fel}). The question then arises as to whether one is truly studying the history of old age itself or the history of its representations -- two things that do not necessarily coincide. Yet, for our purposes -- namely, developing a theory of old age--this distinction is not fundamental.}. This issue is linked to the advances in hygiene and medicine made over the last century; consequently, one cannot glean much from the philosophies of past centuries regarding a question that is, in reality, eminently contemporary.

\section{Physiological and mathematical aspects}

It is in the nature of living beings for biological and cognitive capacities to decline with age -- a fact acknowledged even by the Ancients, including Plato. For instance, human physical performance peaks between the ages of 20 and 30 before entering a decline -- initially gradual, then accelerating -- a process typically modeled by a function such as $f(t) = Ae^{-kt^2}$ or a Gompertz function measuring the exponential decay of the survival rate. Meanwhile, the increase in mortality risk is generally measured using the Gompertz-Makeham law ($\mu _{x}=A+Bc^{x}$); this law posits that the mortality rate is the sum of age-independent terms ($A$, or Makeham terms) and age-dependent terms (the Gompertz function: $\mu _{x}=Bc^{x}$) -- a formulation that represents the inverse of a declining function rather than the decline itself. In contrast, attributes such as memory or information-processing speed are sometimes modeled as undergoing a linear decline followed by a rapid drop-off. Developing a mathematics of aging would, however, require completing three preliminary steps: 1) precisely defining what is to be measured (strength, endurance, neuronal plasticity, etc.) and how these quantities decline; 2) selecting a parametric family of rapidly decaying functions (exponential, Gaussian, Lorentzian, or truncated power-law decay); and finally, 3) verifying against real-world data whether these models fit reality better than classical approaches (such as Gompertz, Weibull, or frailty models). However, these latter laws remain indirect measures of aging\footnote{For instance, the Gompertz law has been used to model tumor growth--cell populations developing in a confined space with limited nutrient availability. The Weibull law (actually a family of distributions) is primarily used in failure analysis for materials. As for frailty models, they are mainly useful in epidemiological contexts where individual variability must be taken into account (see \cite{Ron}).}. A more focused approach would be to employ what are known in mathematics as "declining functions" -- an approach we will adopt here initially. That said, aging is not a purely monotonic decline (there are plateaus, compensatory mechanisms, and vast individual disparities). Furthermore, rapidly decaying functions often imply near-total disappearance within a finite (or asymptotic) timeframe, whereas, biologically speaking, capacities deteriorate but do not all drop to zero at the same age. These aspects must therefore be taken into account.

\section{Schwartz functions}

The precise definition is as follows:

\begin{dfn}
Schwartz functions are defined as functions that are infinitely differentiable and rapidly decreasing, along with all their derivatives of all orders. In French, these functions are known as "fonctions déclinantes" (Declining functions)."
\end{dfn}

Here is an example of a rapidly decreasing function -- specifically, a two-dimensional Gaussian function.

\begin{figure}[h] %  figure placement: here, top, bottom, or page
\hspace{7\baselineskip}
\vspace{0\baselineskip}
\includegraphics[width=4in]{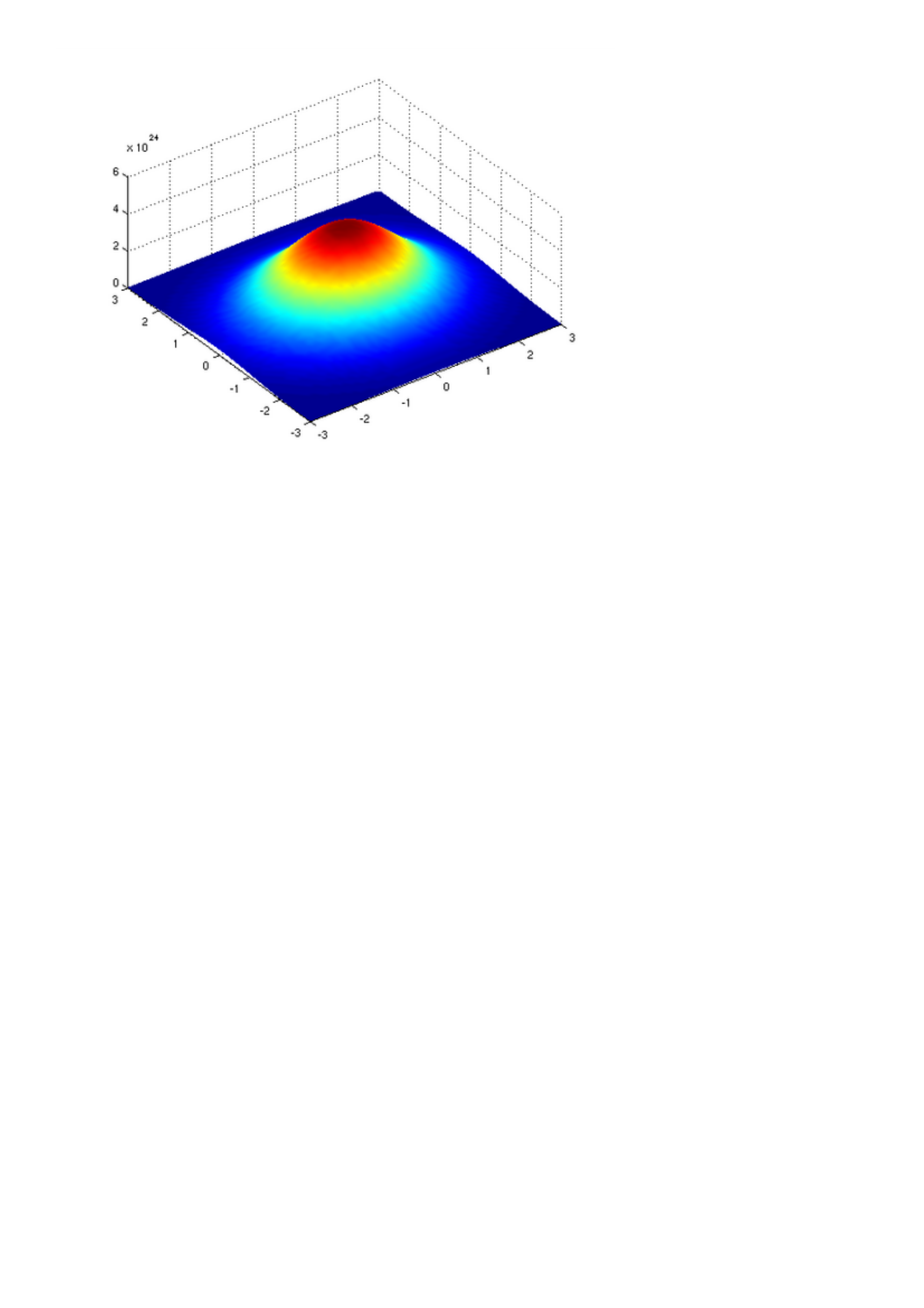}
\vspace{-16\baselineskip}
\caption{Two-dimensional Gaussian function}
\label{fig: imm1}
\end{figure}

Let us now introduce the concept of the "Schwartz space."

\begin{dfn}
The "Schwartz space," denoted by $\mathcal{S}$, is the space of rapidly decreasing functions (i.e., as previously stated, functions that are infinitely differentiable and rapidly decreasing, along with all their derivatives of all orders). The dual space $\mathcal{S}^{\prime}$ of this space is the space of tempered distributions.
\end{dfn}

It is well known that the spaces $\mathcal{S}$ and $\mathcal{S}^{\prime}$ play a fundamental role in the theory of the Fourier transform (see \cite{Sch}). Let us recall one possible definition of the Fourier transform:

\begin{dfn}
The Fourier transform $\mathcal{F}$ is an operation that transforms a function integrable on $\mathbb{R}$ into another function, describing the latter's frequency spectrum. If $f$ is an integrable function on $\mathbb{R}$, its Fourier transform is the function:
\[
F(f) = \hat {f}
\]
given by the formula:
\[
\mathcal {F}(f):\xi \mapsto \hat {f}(\xi )=\int _{-\infty }^{+\infty }f(x)\,\mathrm {e} ^{-{\rm {i}}\xi x}\,\mathrm {d} x
\]
Where:

- $\mathcal {F}(f)= \hat {f}$ is the Fourier transform of $f$;

- $f(x)$ is an integrable function on $\mathbb{R}$;

- $\mathrm{i}$ is the imaginary unit;

- $\xi$ is the argument of the Fourier transform.
\end{dfn}

Geometrically, this transformation can be visualized as follows. Let us take the function $F = a_n \cos nx + b_n \sin nx$ as an example. Here are the initial function and its Fourier transform (see Fig. 2):

\begin{figure}[h] %  figure placement: here, top, bottom, or page
	   \hspace{5\baselineskip}
	      \vspace{0\baselineskip}
	   \includegraphics[width=4in]{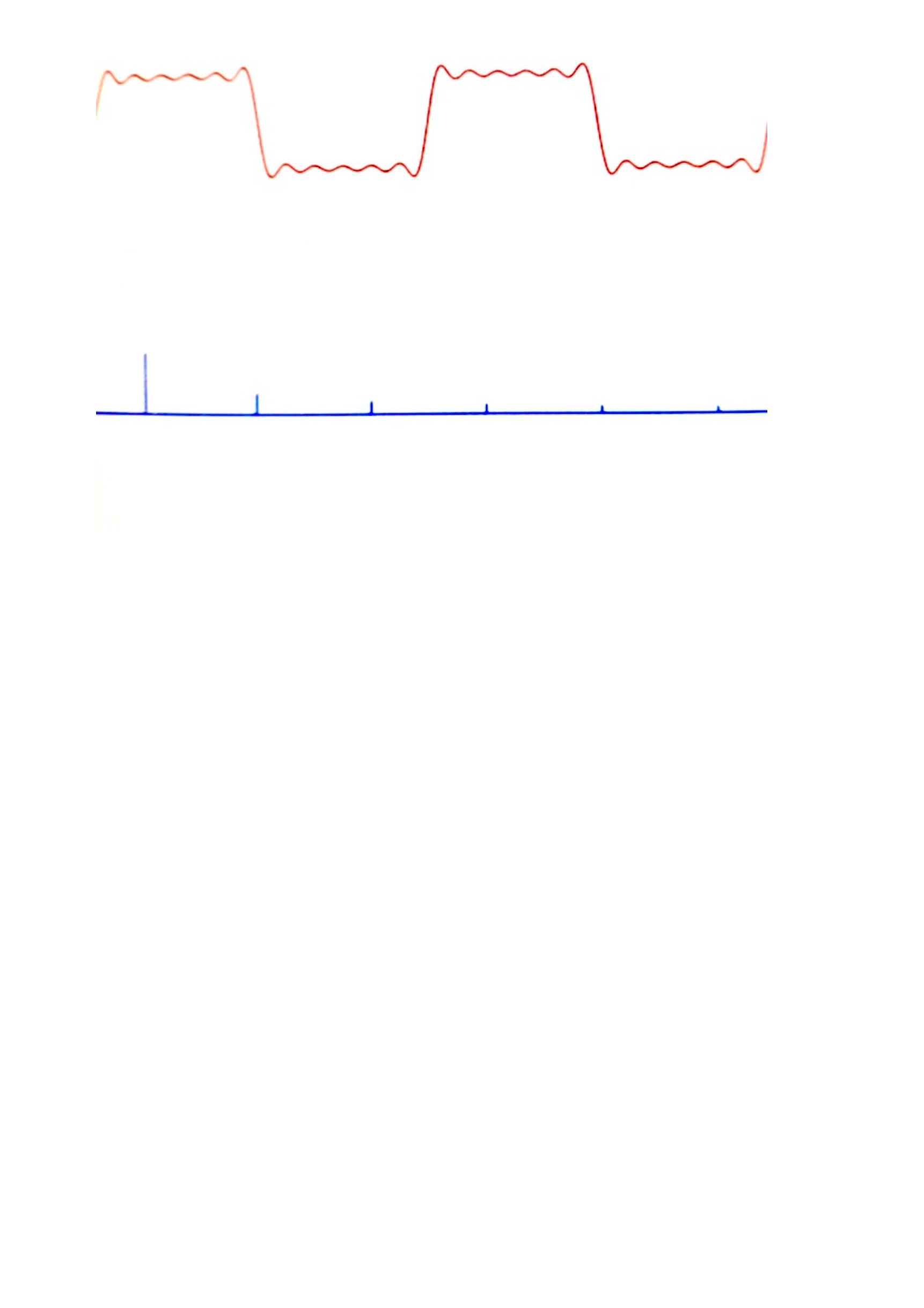}
	      \vspace{-16\baselineskip}
	   \caption{The initial function and its Fourier transform}
	   \label{fig: imm2}
	\end{figure}

And here is exactly how one moves from the function to its transform (see Fig. 3):

\begin{figure}[h] %  figure placement: here, top, bottom, or page
	      \vspace{-1\baselineskip}
	       \hspace{6\baselineskip}
	   \includegraphics[width=4in]{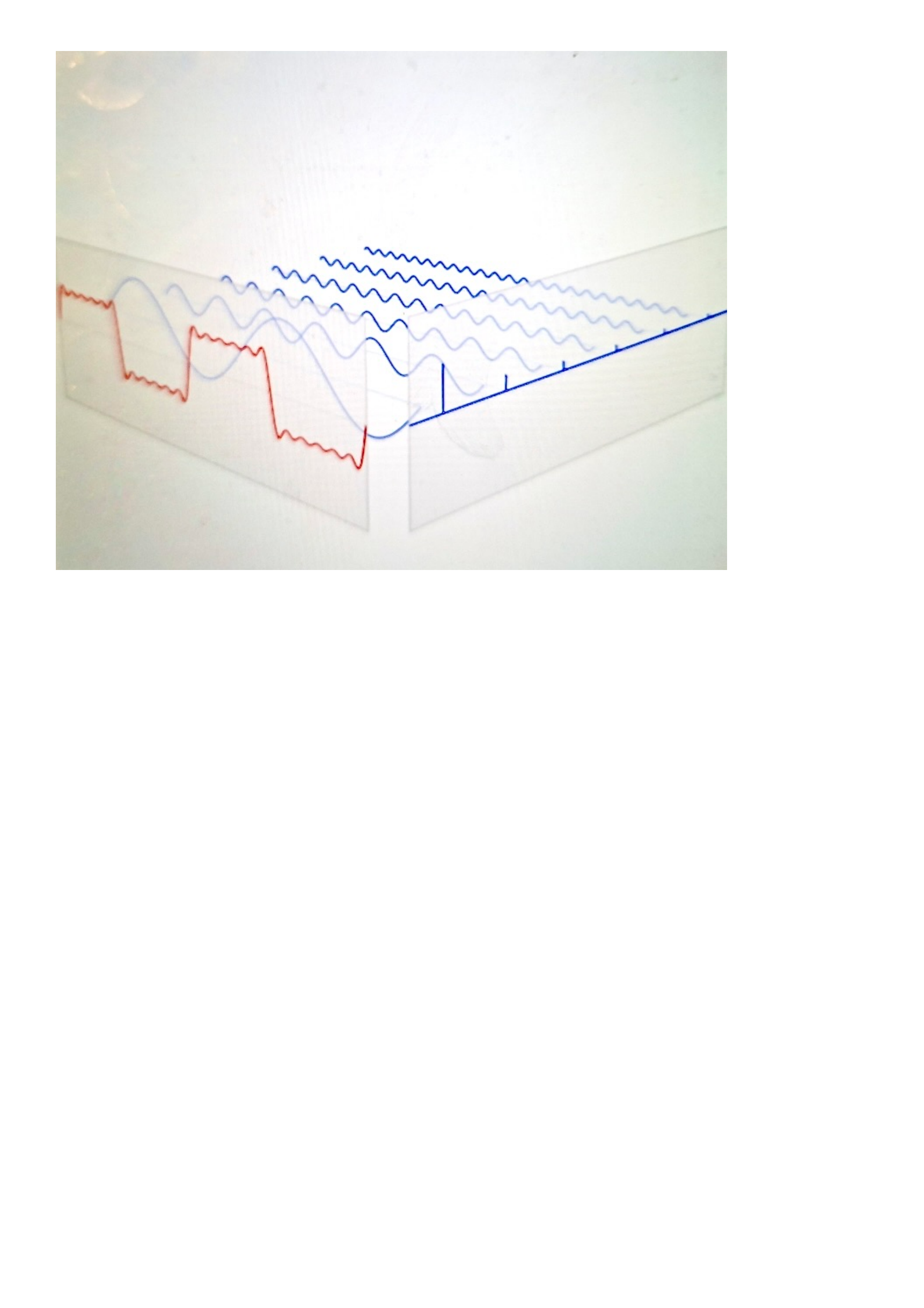}
	      \vspace{-14\baselineskip}
	   \caption{Effect of the Fourier transform}
	   \label{fig: imm3}
	\end{figure}

	It is clear that, when defined on the Schwartz space of $C^\infty$ functions with integrable derivatives, the Fourier transform is a rapidly decreasing function that vanishes at infinity. Conversely, since every Schwartz function is rapidly decreasing, its Fourier transform is of class $C^\infty$. Plancherel's theorem generalizes this situation to square-integrable ($L^2$) functions, allowing us to show that the Fourier transform induces an automorphism of the Schwartz space $\mathcal{S}$ (a vector subspace of $L^2$) onto itself.

The utility of such a formalism for our purposes is evident. If certain decreasing functions associated with aging could be identified as Fourier transforms, it would be interesting to determine what concrete measures might be associated with the inverse Fourier transform, thereby allowing the recovery of the initial function from its spectrum. We will begin by attempting to develop this analogy between mathematical and physiological decreasing functions\footnote{However, it is difficult to imagine that this could lead to any process of individual or social regeneration, or to a palingenesis -- the dream of Orphism, Christianity, and the utopian socialists (Ballanche, Leroux, Reynaud, Chenavard, etc.): the inverse Fourier transform creates nothing. In physics, when moving from a spectrum to a wave, the wave must already exist.}.

\section{Analytical approach}

Suppose we wish to create a
mathematics of aging
based on the theory of "decreasing functions"--functions that vanish at infinity and decrease rapidly. How should we proceed?

The idea is potentially interesting because
aging is precisely a phenomenon in which
certain biological quantities (functional reserve, repair capacity, cellular plasticity, etc.) tend to diminish over time.

Decreasing functions thus offer a
natural language for constructing a mathematical theory of aging. Several possible levels can be distinguished.

\subsection{Aging as the decline of a global capacity}

\subsubsection{The concept of "vitality"}

The idea of ​​a global capacity to maintain oneself in a state of life and good health is an old one.

Historically, the Western concept of "vitality" certainly finds its origins in ancient medical traditions such as Chinese and Ayurvedic medicine, where the balance of vital energies (for example, {\it Qi}\footnote{{\it Qi} (Simplified Chinese: 气; Traditional Chinese: 氣; pinyin: *qì*), a term difficult to translate, is a "natural energy flow" denoting, in Chinese culture, a fundamental principle that shapes and animates the universe and life itself. The difficulty in translation stems partly from the fact that its meaning has evolved over time, influenced by various schools of thought. A quick analysis of the character 氣 (unsimplified script) reveals the symbol for steam (气) above the symbol for rice (米); this yields a highly reductive etymological translation -- "energy produced by the consumption of rice" -- conveying the idea that {\it Qi}  is generated by air and food. However, nutrition is merely one of several ways to produce {\it Qi}. Modern Chinese has retained only the upper component (气), thereby aligning conceptually with the primitive character formed by three horizontal lines (三) -- symbolizing atmospheric currents -- which resembles the character 三 used for the number "three." The concept of {\it Qi} operates on three levels: that of living beings, that of the universe's structure, and that of spirituality. By extension, it is also used to describe a sense of harmony, whether artistic, architectural, or bodily. Interpreting {\it Qi} in terms of "energy" is a specifically Western approach; Chinese texts, conversely, tend to characterize it as a breath or an essence.} or {\it Prana}\footnote{{\it Prana} is a Sanskrit term. As with the Chinese {\it Qi}, the meaning of this compound term is complex, simultaneously integrating the concepts of breath and vital principle into a single synthetic notion expressed by one word. According to the Upanishads, {\it Prana} is a universal vital energy that permeates everything and is absorbed by living beings through the air they breathe.} is considered central to health and well-being but must also, according to these traditions, be cultivated through specific practices such as yoga, meditation, {\it Qi gong}\footnote{{\it Qi gong} (Simplified Chinese: 气功; Traditional Chinese: 氣功; pinyin: {\it qìgong}) is a traditional Chinese form of exercise and breathwork based on the knowledge and mastery of breath, combining slow movements, breathing exercises, and concentration. The term literally means "achievement or accomplishment ({\it gong}) related to {\it Qi}" or "mastery of breath,".} etc. In this context, vitality is viewed as an essential force flowing through body and mind, and is inherent to health (see \cite{Bra}).
Along similar lines, the Dutch philosopher Baruch Spinoza regarded the notion of {\it conatus} -- or the effort to persevere in one's own being -- as a characteristic of all individuality\footnote{Spinoza, {\it Ethics}, III, Proposition 6: "{\it Unaquaeque res, quantum in se est, in suo esse perseverat conatur}" (each thing, insofar as it lies within itself, strives to persevere in its being), and Proposition 7: "{\it Conatus, quo unaquaeque res in suo esse perseverare conatur, nihil est praeter ipsius rei actualem essentiam}" (the effort by which each thing strives to persevere in its being is nothing other than the actual essence of that thing).}, though this notion was physical rather than specifically vital. Schopenhauer's "will-to-live" or, {\it a fortiori}, The Nietzschean "will to power" -- of which these are avatars -- certainly intensified this force, rendering it, in the process, more positive or, let us say, more active. Henri Bergson, too, would inherit some of its meanings through the concept of the "élan vital" -- central in {\it Creative Evolution} (see \cite{Ber}, 569–578, 708, 713–714), though its sense is oriented more toward the evolution of species than of individuals -- while adding to it, as a famous passage from {\it The Two Sources of Morality and Religion} demonstrates, the notions of unpredictability, an internally felt indivisibility, and a real, efficacious duration (see \cite{Ber}, 1072). 

Fitzgerald, for his part, would highlight the unequal distribution of vitality among men and the possibility of its sudden failure.
For him, vitality -- and this is the great tragedy of "The Crack-up" -- is what cannot be communicated: "Vitality never {\it takes}," he writes. "You either have it or you don't, like you have health, or brown eyes, or honor, or a baritone voice" (see \cite{Fit}, 483)\footnote{[O]f all natural forces, vitality is the incommunicable one. In days when juice came into one as an article without duty, one tried to distribute it -- but always without success; to further mix metaphors, vitality never “takes”. You have it or you haven't it, like health or brown eyes or honor or a baritone voice (*The Crack-up*).}. Furthermore, against a naturalism that reduces the vitality of living beings to a positive and measurable constitution, Georges Canguilhem, just like the German Helmuth Plessner\footnote{Helmuth Plessner (1892-1985), German philosopher and sociologist, is, with Max Scheler and Arnold Gehlen, one of the main representatives of philosophical anthropology. He is notably the author of various works, including *Die Stufen des Organischen und der Mensch*, an introduction to his anthropological philosophy, dating from 1928.}. Although in different terms, he posited the irreducibility of this vitality to the positive discourse of science (see \cite{Ebk})\footnote{Undoubtedly, unlike Fitzgerald, who treats vitality more as a kind of intrinsic force, which tends to essentialize it, Canguilhem considers vital action primarily as a "debate" between the living being and its environment, health, always fragile, being above all the capacity to become ill and to recover (see \cite{Can}, 187).}.
Although somewhat vague, if not obscure, the notion of vitality is therefore a frequently used term, and even if the existing literature on this notion does not provide a concise and unique definition of what it actually is\footnote{Although close, it seems broader, in particular, than the Freudian notion of libido, which is too heavily tinged with sexuality.}, the experience of an energy, both psychic and physical, felt by a person and available to them seems sufficiently widespread that we can consider retaining it, at least as a first approximation, as the basis for a mathematical representation of certain behaviors of living beings.

It therefore seems possible to introduce a function $V(t)$
representing this "vitality" or, if you prefer, its equivalent: the "physiological reserve" of strength or energy of a living being at age $t$\footnote{This could be assessed via a discrete or continuous scale measuring not only physiological variables (muscle tone, hormonal balance, immune response to viral infections, etc.), but also the subject's cognitive and memory capacities, as well as their ability to overcome adversity or difficult situations (resilience) and, if necessary, to engage in new projects.}.

We impose:

\begin{itemize}
\item $V(t) > 0$,
\item $V'(t) < 0$,
\item $\lim_{t \to \infty} V(t) = 0$

\end{itemize}

For example:

\begin{equation}
V(t) = V_0 e^{-\lambda t}
\end{equation}

This law expresses a loss proportional to what
remains available.

But biologically, actual aging is
not necessarily exponential.

\subsubsection{Rapid decay and collapse
of reserves}

A rapidly decaying function -- in the sense used in mathematical analysis -- decays faster than any power law; that is:
\begin{equation}
\forall n, t^n V(t) \to 0.
\end{equation}
The classic example is:
\begin{equation}
V(t) = e^{-t}
\end{equation}
One interpretation is that youth corresponds to a relatively stable phase, after which residual capacities become extremely low once a certain threshold is crossed. This would yield a sort of theory of "evanescent biological reserves."

\subsubsection{Old age as the dynamics of
multiple declining functions}

A single variable, however, is likely too limited to capture the full richness of the concept of "vitality." We could therefore introduce several variables, such as:

\begin{itemize}
\item $M(t)$: memory,
\item $I(t)$: immunity,
\item $R(t)$: repair capacity,
\item $S(t)$: muscle strength.
\end{itemize}

Each function would then decline at its own rate while remaining positive until the moment of death, at which point it would -- naturally -- drop to zero:
\begin{equation}
M(t), I(t), R(t), S(t) > 0.
\end{equation}

Old age would then become a trajectory within a functional space, allowing us to define a biological age $A(t)$ as a function of time, such that:
\begin{equation}
A(t) = F(M (t), I(t), R(t), S(t)).
\end{equation}

This brings us closer to the concept of a differential geometry of life, or to the theories of evolutionary dynamical systems for which such geometry serves as the underlying formalization.

\subsubsection{Using Schwartz spaces}

However, explicitly employing rapidly decaying functions suggests going a step further and making full use of the previously introduced concept of the {\it Schwartz space}:

\begin{dfn}
The Schwartz space consists of functions such that:
\begin{equation}
sup_x |x^\alpha D^\beta f(x))| < \infty.
\end{equation}
\end{dfn}

This can be interpreted as follows:

\begin{itemize}
\item biological mechanisms possess frequency components;
\item aging acts as a filter that progressively eliminates the high frequencies of physiological activity.
\end{itemize}

In other words, a young organism exhibits complex fluctuations across multiple time scales, whereas an aging organism becomes more rigid and less variable.

This idea is, in fact, quite similar to certain observations in the physiology of heart rhythm and biological variability.

\subsubsection{An entropy-inspired approach}

Another approach to aging is also possible. Instead of directly modeling what declines, one could -- conversely -- model what increases by setting:
\begin{equation}
D(t)=1-V(t),
\end{equation}
where $D$ represents the accumulation of damage.

Aging would then correspond to the monotonic increase of $D$ (comparable in this respect to the increase of entropy in a closed system), while $V$ would be its decreasing transform. This reveals a duality between a vitality tending toward 0 and a set of damages (or increasing disorders) tending toward 1.

\subsubsection{A more abstract formulation}

In a more abstract formulation, the most novel idea would be to
define aging not as an "age," but
as an asymptotic equivalence class
of functions.

Two organisms would share the same type of
aging if their vital functions satisfy:

\begin{equation}
\frac{V_1(t)}{V_2(t)} \to C>0.
\end{equation}

They would then belong to what might be called the same "asymptotic
aging species."

One could thus classify these forms of aging, distinguishing, for example, between:

\begin{itemize}
\item polynomial aging: $t^{-\alpha}$;
\item exponential aging: $e^{-\lambda t}$;
\item super-exponential aging: $e^{-t^2}$;
\item stretched aging: $e^{-t^\beta}$.
\end{itemize}

This would yield a genuine "analytical
theory of forms of senescence."

\section{Algebraic approach}

That said, if one were truly to construct an authentic "mathematics of
aging," it might be preferable -- while initially remaining at the quantitative level -- to start not from the simple notion
of a decreasing function, but from the more appropriate concept of a
{\it dissipative semigroup}: the physiological state
$X(t)$ would then evolve under the action of an operator that
progressively destroys biological
information.

Declining functions that vanish at
infinity or exhibit rapid decay would thus appear
as the {\it observables} of this dynamical
system rather than as the fundamental {\it objects}
themselves.

This would make it possible
to simultaneously integrate the loss of function,
the accumulation of damage, and the interactions
between organs within a unified framework.

\subsection{"Inclines"}

For about forty years, structures
known as "inclines" have been recognized;
these could be linked to
this new line of inquiry. An "incline" is, literally, a set $\mathcal{K}$
on which two binary operations, $+$ and $\cdot$, are defined,
satisfying the following axioms ($a, b, c \in \mathcal{K}$):

\begin{enumerate}
\item $+$ is commutative: $a+b = b+a$;
\item $+$ and $\cdot$ are associative: $a+(b+c) = (a+b)+c, a(bc) = (ab)c;$
\item $\cdot$ is distributive over $+$: $a(b+c) = ab + ac, (b+c)a = ba+ca$;
\item $+$ is idempotent: $a+a = a$;
\item the specific property of "inclines" holds: $a+ac=a, c+ac=c$. \end{enumerate}

These "inclines" -- idempotent algebras introduced by Cao, Kim, and Roush
in the 1980s (see \cite{Cao}, 1–2) -- are structures
closely related to semirings\footnote{In fact, we have the following inclusions: Boolean algebras $\subset$ fuzzy algebras $\subset$ distributive lattices $\subset$ inclines $\subset$ semirings (see \cite{Cao}, 2).} and in which the product acts
as an operation of diminution or
contraction.

One of their characteristic properties is that this
product essentially satisfies:
\begin{equation}
ab \le a,  \qquad  ab \le b.
\end{equation}
In other words, combining two states can never yield something greater than either of the two initial states, and this is precisely what evokes the concept of aging.

In a theory of aging, one could interpret the elements of the "incline" as levels of biological capacity, the ordering as an ordering of vitality, and the product as an interaction between physiological functions.

In this context, property (9) would mean that a physiological interaction never spontaneously creates additional reserve; at best, it merely preserves a portion of what already exists.

This differs significantly from classical algebra, where interactions can amplify.

However, there is a deeper element to these structures.

Indeed, "inclines," as previously noted, are characterized by idempotent operations of the type:
\begin{equation}
a+ab = a,
\end{equation}
where multiplication is compatible with a natural order\footnote{The concept of an "incline," by incorporating this notion of order, would naturally allow one to speak of what is "older" or "more degraded." However, it also enables the description of degradation not merely through a single value, but through a network of dynamic relationships between a system's components (see, for example, \cite{Cao2}).}.
In this sense, these algebras are closely related to idempotent semirings (see \cite{Gon}) and certain constructions in tropical geometry (see \cite{Pla})\footnote{Recall that tropical geometry (the term "tropical" merely alludes to Brazil, the home country of mathematician Imre Simon, a pioneer in the field) refers to the possibility of defining two algebras on the set of real numbers -- the min-plus algebra and the max-plus algebra. Both possess the structure of a commutative semifield and are isomorphic to one another via the mapping $x\mapsto -x$. The choice of the reference algebra depends on the authors and the specific field of study.
For the max-plus algebra, the definition is as follows:
$a\oplus b=\max(a,b)$. (For the min-plus algebra, the maximum is replaced by the minimum).
Similarly, tropical multiplication (or the tropical product) $\odot$ (or $\otimes$) is defined as:
$a \odot b=a+b$.
Thus, the result of the tropical multiplication of two numbers is their standard sum.
Like standard addition, tropical addition is commutative and associative. However, there is no neutral element in $\mathbb{R}$ (to have one, one would need to work in $\mathbb{R} \cup \{+\infty\}$).
Tropical multiplication is, like standard multiplication, commutative and associative. It is also distributive with respect to tropical addition.
The number 0 is the neutral element for this multiplication, and every element has an inverse under this operation.
Since the first operation lacks a neutral element and the existence of symmetric elements, the structure as a whole is not a field. It is referred to simply as the "semifield" $(\mathbb{R}, \oplus, \odot)$.}.

Property (9) thus aligns with a logic of saturation rather than accumulation.
Biologically speaking, it evokes phenomena where:

\begin{itemize}
\item additional repair yields no further benefit beyond a certain threshold;
\item functional reserve is bounded;
\item certain mechanisms are limited by physiological ceilings.
\end{itemize}

In most standard physical theories, one works with groups, rings, and vector spaces, where the defined operations allow for both increase and decrease. "Inclines," by contrast, introduce an algebra intrinsically oriented toward loss.

This is precisely the type of structure one would seek if aiming to formalize senescence, the loss of biological information, or the reduction of physiological degrees of freedom. Thus, one interpretation of property (10) could be that time ($b$), when applied to a capacity ($a$), cannot increase it.

\subsection{The concept of algebraic age}
In this context, one could introduce an *algebraic age* defined via iteration.
Indeed, suppose that at each age $t$ there exists a set of accessible states $A_t$.
Aging could be described by the equation:
\begin{equation}
A_{t+1} \subseteq  A_t.
\end{equation}
Alternatively -- or as another solution -- if we consider a partially ordered set of states $\mathcal{L}$ equipped with an ordering relation $\le$, and if:
\begin{equation}
X(t) \in \mathcal{L}
\end{equation}
represents the state of the organism at time $t$, we might assume that:
\begin{equation}
X(t+1) < X(t). 
\end{equation}
In other words, aging would constitute a descending sequence within $\mathcal{L}$.

Algebraic age could then be defined as the {\it ordinal distance} from a maximal reference state $X_0$.
The underlying idea is that available biological information progressively contracts.
In this view, old age is less a matter of functions tending toward zero than a matter of the monotonic reduction of a range of possibilities.

\section{Dissipative systems and algebra}

A dissipative system (or dissipative structure) is a system that evolves within an environment with which it exchanges energy or matter. It is, therefore, an open system operating far from thermodynamic equilibrium. A dissipative system is characterized by the balance of its exchanges (energy exchange, entropy production) and the spontaneous emergence of spatial symmetry breaking (anisotropy), which can sometimes give rise to a complex, chaotic structure. A limit to this degradation lies in the fact that the system can be momentarily stabilized through its "consumption" of energy derived from the environment. The term "dissipative structures," coined by the Russian-born Belgian physicist and chemist Ilya Prigogine in the 1970s (see \cite{Pri1}; \cite{Pri2}; \cite{Pri3}), has made it possible to highlight -- far from equilibrium -- the emergence of novel structures that underpin many forms found in nature (see \cite{Par}). A simple example is that of Bénard cells. More complex examples include reactions such as the Beloussov-Zhabotinsky reaction and, of course, life itself; this justifies the mention of dissipative structures here and may -- later on (see section 9) -- allow us to derive a different conception of aging based on them.

Since, within an "incline," the powers tend to become increasingly small in the associated order, the sequence of iterates naturally represents a trajectory of decline.
Here, we encounter the concepts of decreasing functions and objects that vanish asymptotically.

There are, therefore, at least two ways to mathematize aging:

-- via Analysis: decreasing functions, asymptotic limits;

-- via Order Algebra: "inclines," semirings, idempotents, tropical algebras.

The first describes the temporal form of decline -- the asymptotic disappearance of capabilities.
The second describes the structural rules that render this decline inevitable.
"Inclines," in particular, perfectly describe the contraction of the structure of possible states.
It is perhaps this second approach that is, philosophically, the most novel: aging is no longer viewed as a specific function, but rather as a fundamental algebraic property of the system itself. One might be tempted to elaborate on this in greater detail.

\section{Old age as an algebraic property}

Let us return to "inclines" -- an algebra of idempotent semirings.
There is undoubtedly a possible conceptual link
to dissipative structures, even if it is not immediately obvious.
As previously mentioned, "inclines" are ordered structures in which
addition is idempotent and
multiplication is compatible with a natural order.
This is highly significant.
Indeed, regarding the question of aging, what appears interesting
is not so much the technical definition of "inclines" as their underlying philosophy.
Note that:

1. order is fundamental to them;

2. operations tend to select states rather than
add them together in the classical sense;

3. many phenomena are described in
terms of dominance, absorption, and
loss of information. Yet aging can be viewed precisely
as a dynamic process in which:
\begin{equation}
X_{t_1} < X_t,
\end{equation}
that is, a sequence of increasingly constrained
physiological states.

As previously discussed, one could even imagine that the state of an
organism is not a real number but an
element of a (partially) ordered set $\mathcal{L}$, representing
the set of possible biological capacities.
Old age would then no longer be a
numerical decline but, as suggested earlier, a "descent"
within this partial order.

However, another point appears even more
profound.

Classical decreasing functions
describe a quantitative loss (where $V(t) \to 0$).
"Inclines," by contrast, suggest a loss that is
{\it structural} in nature. In other words: in an analytical approach,
one loses amplitude, whereas in an "incline"-based approach,
one loses possibilities.

An aging organism is, in fact, not merely a
"weaker" organism; it is an organism
whose space of accessible states is shrinking.
Old age closely resembles a dynamic
collapse of physiological variety.
Yet, in this view, it is less a
matter of functions tending toward zero than a
matter of the monotonic reduction of a range of
possibilities.

One could therefore argue that rapidly decreasing
functions and "inclines" offer perspectives on the
same phenomenon, but from two different angles:

\begin{itemize}
\item functions describe the asymptotic
disappearance of capacities;
\item "inclines" describe the contraction of the
structure of possible states.
\end{itemize}

The first approach is analytical; the second
is ordinal.

That said, if the goal is to construct a genuine "mathematics
of aging," one might suspect that the second approach
could be more fundamental than the first: what
characterizes an aging organism
is perhaps not so much the decline in overall vitality values ​​as the
progressive reduction in the richness of the
configurations it can still realize.

\section{Clarifications on the concept of Age}

Here, we shall attempt to reformulate, with even greater precision,
the idea naturally emerging from the
preceding reflections on ordered structures and
"inclines."

\subsection{Chronological age versus algebraic
age}

Chronological age is simply time $t$. Algebraic age can be viewed as a measure of an organism's
position within an ordering structure
representing its capabilities.
Let us once again consider the set $\mathcal{L}$ of physiological
states, equipped with an ordering relation.

If $X(t) \in \mathcal{L}$ is the organism's state at time
$t$, we assume that:
\begin{equation}
X(t+1) < X(t).
\end{equation}
In other words, aging appears as a descending sequence
in $\mathcal{L}$.

As previously mentioned, algebraic age could then be defined as
the ordinal distance to a maximal reference
state $X_0$. However, this is merely one initial possibility. Let us now introduce other, equally valid definitions.

\subsection{Definition based on height}

Assume that $\mathcal{L}$ is a lattice. If this lattice $\mathcal{L}$ is discrete, we can define:

$A(X)$ = maximum length of a descending chain starting from $X$.

For example:
\begin{equation}
X_0> X_1>X_2>X_3> 0.
\end{equation}

Then:
\begin{equation}
A(X_0) = 0, \qquad A(X_1) = 1, \qquad A(X_2) = 2, \qquad A(X_3) = 3.
\end{equation}
With this approach, age is no longer measured in years but in
terms of structural losses.

\subsection{Definition based on rank}

In a more quantitative version, one would associate a
rank function with each state:
\begin{equation}
r \to \mathbb{R}^+,
\end{equation}
where:
\begin{equation}
X \le Y \Rightarrow r(X) \le r(Y). \end{equation}
We then define:
\begin{equation}
A(X) =r(X_0) - r(X).
\end{equation}
Algebraic age would thus measure the amount of structure lost since the initial state.

\subsection{Introduction of an operator $T$}
With the iteration suggested above, the idea was to replace a continuous aging function
with an operator.
We thus introduce an operator:
\begin{equation}
T: \mathcal{L} \to \mathcal{L},
\end{equation}
such that:
\begin{equation}
T(X) \le X.
\end{equation}

A unit of time corresponds to the application of $T$.
The evolution then becomes:
\begin{equation}
X_{n+1} = T(X_n).
\end{equation}
or alternatively
\begin{equation}
X_n = T^n (X_0).
\end{equation}

In this version, chronological age is $n$, but algebraic age
is determined by the position of $T^n(X_0)$ within
the ordering.

\section{A possible synthesis}

Returning to our initial idea, we can
take a vitality function $V$ and define:
\begin{equation}
T(V)(t) = \phi(t)V(t),
\end{equation}
where:
\begin{equation}
0 \le \phi(t) \le 1.
\end{equation}

Each iteration then applies an additional
loss, such that after $n$ steps, we have:
\begin{equation}
T'(V)(t) = \phi(t)^n V(t).
\end{equation}
Algebraic age could then be redefined as:
\begin{equation}
A(V) = -\log\left(\frac{\|V\|}{\|V_0\|}\right),
\end{equation}
where $\|\cdot\|$ is an appropriate norm.
This yields an additive quantity:
\begin{equation}
A(T^n V) = A(V) + nc,
\end{equation}
in the case where the operator always contracts
by the same factor.

However, arguably the most elegant version of the intended formalization would be based on the following observation.
If one wished to truly link old age, order, and
"inclines," one could define algebraic age as:
\begin{equation}
A(X) = \min\{n : T^n(X_0) \le X\}.
\end{equation}
In other words, the age of an organism's state would be the minimum number
of iterations of the degradation process
required to reach or surpass that
state within the ordering.

Note that this definition relies neither on a metric nor
on a notion of distance. It relies
solely on the ordering and iteration of a
monotone decreasing operator, making it
very close in spirit to "inclines" and
ordered algebraic structures.

\section{Algebraic age and the age of a relational structure}

For anyone familiar with graph theory and the theory of relations, the concept of {\it algebraic age} inevitably brings to mind
the concept of the {\it age of a
relational structure} in the sense of
Cameron (see \cite{Cam}).

This connection is, in fact, more
than a mere terminological analogy.
In the theory of relational structures,
the "age" of a structure -- as defined by Peter J.
Cameron (see \cite{Cam}, p. 50), perhaps more clearly than by Fraïssé (see \cite{Fra}, p. 279) -- is the set of finite structures
that appear as induced substructures
of a given structure. It is thus a
concept of combinatorial heritage: the age describes
everything that can still be realized locally
within the structure. Formally, we have the following definitions:

\begin{dfn}[Relational structure]
A relational structure $(X, R)$ consists of a set $X$ and a family $R = \{r_\alpha: \alpha \in A\}$ of finitary relations on $X$. In other words, $r_\alpha$ is a subset of $X^n$, where $n = n(\alpha) \ge 1$ is the arity of $r_\alpha$.
\end{dfn}

From this point, Cameron introduces the notion of a "substructure" of a relational structure.

\begin{dfn}[Substructure]
A substructure of $(X, R)$ is a relational structure $(Y, S)$ where $Y$ is a subset of $X$ and $S = \{s_\alpha: \alpha \in A\}$, with $s_\alpha = r_\alpha \cap Y^n$ and $n = n(\alpha)$.
\end{dfn}

The domain of a relational structure in Cameron's sense is always finite or countable. The author subsequently introduces the notion of the "age of a structure" as follows:

\begin{dfn}[Age of a structure]
The {\it age} of a relational structure $X$, denoted $Age(X)$, is the class of all finite structures that are isomorphic to substructures of $X$.
\end{dfn}

In this context, the age is, of course, a proper class rather than a set, which causes no particular difficulty.

If we transpose this to the case of biological aging,
we could say:

\begin{itemize}
\item at a given time $t$, the organism possesses a relational structure $R_t$;
\item its "age," in Cameron's sense, is the set
of physiological configurations that are still realizable.
\end{itemize}

The difference between the two definitions of age is subtle but important.
For Cameron, $Age(R)$
is a set of finite structures.
In the previous outline, algebraic age
was a number or a rank, denoted $A(X)$.

However, one can obviously merge the two ideas: instead of defining an organism's age by a
quantity, It is defined by its class of
configurations that are still present.
A young organism possesses a very rich
combinatorial age:

\begin{itemize}
\item many possible immune responses;
\item many accessible neuronal patterns;
\item many modes of adaptation.
\end{itemize}

With aging, certain configurations
disappear.
Age then becomes an ordered entity.

Thus, a young heart can produce a wide
variety of adaptive rhythms, whereas an
aged heart often exhibits
more rigid dynamics. In the language defined here, one could say that the
relational age of the cardiac system decreases.

There is even a construct that seems
particularly promising.

Let us assume that the physiological state is described
by a large graph of the type $G_t$. The vertices represent cells,
tissues, or functions; the edges represent
interactions. We then define $Age(G_t)$
as the Cameron age of the graph.

Aging would thus be observed through the
gradual disappearance of certain types of
subgraphs, and we would have:
\begin{equation}
Age(G_{t+1}) \subseteq Age(G_t).
\end{equation}
Biological age would literally become a
loss of relational patterns.

In summary, one could say that, in classical aging theory, one seeks a function $V(t)$ that decreases,
whereas in an algebraic and structural approach to age in the Cameron sense, one
would instead seek a decreasing family of
relational classes of the type $A_t$. Time would no longer measure how much
vitality remains, but rather how many structures are still
realizable. This is why invoking this algebraic age
sheds more light on the problem than the
rapidly decaying functions themselves.

A function $V(t) \to 0$ describes a loss of
{\it quantity}. An age in the Cameron sense describes a
loss of {\it structural complexity}. Yet
actual aging often resembles the latter
more than the former. One does not merely become
"less capable"; one becomes
capable of fewer and fewer types of
behaviors, repairs, and adaptations. This is precisely the kind of phenomenon that relational-structure ages were designed to capture, even in a completely different context.

\section{Another perspective on old age}

From Cicero (see \cite{Cic}) to Hermann Hesse (see \cite{Hes}), people have often tried to ward off old age -- with its declining functions and the general weakening of vitality we have described. 

In Cicero’s work, it is the elderly Cato who speaks: of the four flaws commonly attributed to old age (unfitness for affairs, the gradual wearing down of strength, the total absence of pleasure, and finally, the increasingly encroaching proximity of death), he sets out to turn the first into an advantage. 

First, there are affairs suited to old age and within its purview -- and these are perhaps the most important: long experience of the past teaches one to understand the future and fosters a kind of prophetic art -- precisely the province of the elderly, for one must have lived a long time to have seen a great deal. They can thus enlighten the young, mold men, and thereby render a tremendous service to society. 

Next, it is said that in old age, one's strength is worn down. But what kind of strength are we talking about? Physical strength is not truly necessary for the roles reserved for the elderly, who can, moreover, retain a measure of vigor if they have led a temperate and disciplined life. Come what may, of course, one must submit to the laws of nature; indeed, every age has its limits. That said, while the body may be weak in old age, the soul -- that is, the noblest part of a human being -- retains sufficient strength to act, and to act rationally.
	
	Admittedly, dulled senses are no longer suited to savoring certain pleasures. Yet, this is a fortunate thing, for sensual indulgence can prove destructive; a man in the grip of insatiable delights and the most ruinous passions loses the habit of thought. Indeed, no pleasure remains for the man who has lived a life of excess. Conversely, the man who has never squandered his health can still experience pleasant sensations and the quiet joys of a country life that reconnects him with the simplicity of nature. 

Of course, old age borders on death, and this frightening proximity brings sadness. Yet, in truth, death can strike at any age, and the man who has lived in accordance with reason need not fear it. Only he who never truly knew how to live will fail to know how to die. The wise man views it merely as a transition from a wearisome life to a happier state. 

From this perspective, old age calls for a model that goes beyond a mere decline in vital functions or a simple loss of structural complexity. It is a matter, rather, of replacing certain organizational structures with others -- or, more precisely, of selecting from the web of relational structures linking a living being to its environment those that were previously neglected or undervalued, and instead developing them to serve as the foundation or central axis of this new phase of life. Hermann Hesse's "In Praise of Old Age"\footnote{In reality, the text is titled "Mit der Reife wird Man immer Jünger (Betrachtungen und Gedichte über das Alter, edited by Volker Michels)" (in other words: "With age, one grows younger. Reflections and poems on aging, edited by Volker Michels").} is undoubtedly a more poetic version of {\it De Senectute}, yet it confirms the words of the old Stoic sage.

What is life made of as one ages? Plans dwindle, the horizon narrows, the pace of life slows, and new prospects become rare. It is rather like feeling oneself dying while still alive. One would like to possess Hermann Hesse's Eastern wisdom: his sense of quiet contemplation in the face of the precariousness of all things; his joyful acceptance of the melody of the ephemeral; his love of change and metamorphosis -- a love that made him not fear death, but desire it with curiosity -- and, one might even say, almost with eagerness -- like any transition of state: "tomorrow, the day after, soon, very soon," he writes in (\cite{Hes}, 8–9), "I shall be something else; I shall be the foliage, I shall be the earth, I shall be the root; I shall no longer set words down on paper, I shall no longer breathe the scent of the magnificent yellow wallflower, I shall no longer carry a dentist's bill in my pocket, I shall no longer be tormented by formidable bureaucrats regarding my certificate of nationality; I shall be the cloud floating in the air, I shall be the ripple in the stream, I shall be the leaf budding on the shrub; I shall enter into oblivion, I shall plunge into the cycle of metamorphoses so ardently desired."

After fifty years, Hesse tells us, "man ceases to indulge in childish things" (see \cite{Hes}, 30). Only one question seems of interest: that of the Spirit and of faith, of meaning and of the authentic piety which alone enables one to confront suffering and death -- and to rise, so to speak, to the challenge of these events. In this context, old age affords unique experiences, impossible to undergo at any other stage of life.

First and foremost, a mystical experience: "When one is old enough, it feels as though one's entire existence -- with its joys and sufferings, its loves and discoveries, its friendships and affairs, its books, music, travels, and work -- amounts to nothing more than a long detour leading to the unfolding of those moments when God reveals Himself, when the meaning and value of all that exists and comes to pass are disclosed to us through the form of a landscape, a tree, a face, or a flower. It may be that in our youth we admired the spectacle of a blossoming tree, of clouds gathering in the sky, or of a storm with greater passion and fervor; yet, to undergo the experience of which I speak, one must have reached a certain age. One must have seen, lived, thought, felt, and endured an infinite array of things; one must sense the waning of vital instincts, the newfound fragility of the body, and the nearness of death" (see \cite{Hes}, 54-55).

The same observation applies when, watching a poplar swaying in the wind, Hesse comes to experience the harmony of the world. "With joy and without fear -- indeed, with a light heart -- the poplar surrendered its branches and its cloak of leaves to the damp wind, which was gathering considerable strength. The song it sang on that stormy day, the shapes its tapering crown traced against the sky, seemed to me marvelous, incomparable. They expressed joy as well as gravity, active will and submission, the interplay of freedom and destiny" (see \cite{Hes}, 61). And further on: "the ceaselessly shifting dance of the treetop in the storm was but an image revealing the mystery of the world -- beyond strength and weakness, good and evil, acting and suffering. For an instant -- a brief minute of eternity -- I saw revealed in a pure and perfect form (purer and more perfect than if I had read Anaxagoras or Lao-Tzu) that which usually remained hidden and secret. And in that moment, I felt that perceiving this image and deciphering its meaning had required not only the miracle of that spring hour, but also the journeys and wanderings, the follies and experiences, the pleasures and sufferings of decades upon decades" (see \cite{Hes}, 61–62).

Thus, old age -- a source of much pain -- is also the origin of many graces. Declining faculties bring about lethargy and rigidity. Yet, viewed in a different light, "this can also appear as tranquility, patience, humor, a high degree of wisdom -- the Tao" (see \cite{Hes}, 71). Whether this shift in one's relationship with the world -- characteristic of old age -- stems from enlightened wisdom or is a symptom of aging (such as circulatory disorders) remains a mystery, even to the person experiencing it. In any case, "it is only as one ages that one realizes beauty is rare, that one understands the miracle of a flower blooming amidst ruins and cannons, or the survival of literary works amidst newspapers and stock market quotes" (see \cite{Hes}, 72). Thus, old age -- which takes seriously that which transcends the individual -- entails a kind of sacrifice: the elderly person must let go of the self or, as German mystics once put it, "undo their being" (see \cite{Hes}, 72).

Old age thus allows one to attain a certain dignity rooted in accumulated memories and experiences -- a past that persists into the present and helps one overcome the trials of aging.

Another characteristic of old age is a growing fondness for established habits and repetition.

Hesse explains that, unlike young people, the elderly have virtually nothing left to discover. "They have long since lived through the fundamental experiences suited to their nature and destined for them. Consequently, the 'new' -- and increasingly rare -- events in their lives are repetitions of what has already been experienced many times over" (see \cite{Hes}, 139).
Yet repetition encompasses difference: "these events also represent something novel," writes Hesse. "They are certainly not the first of their kind, but they are no less real for that, as each time they transform into a turning inward, into introspection" (see \cite{Hes}, 139–140). Consequently, the unexpected return of a forgotten memory takes on a new force and power -- qualities not possessed by memories that are conscientiously cultivated.

Undeniably -- and unlike youth, which thirsts for the unknown and the new -- "old age develops a taste for established habits and repetition; it constantly seeks to revisit the same places, people, and situations, likely because it aspires to enter the realm of memory and feels an unflagging need to verify what memory has preserved. Perhaps, too, it harbors the wish -- the slender hope -- of seeing this treasure grow, of one day rediscovering the trace of some forgotten and lost event, encounter, image, or face, and adding it to the sum of its memories" (see \cite{Hes}, 141-142).

Finally, a last characteristic: old age no longer seeks truth, nor even love. "Truth is an ideal typical of youth; love is a dream cherished by the mature man -- the one striving to be ready to face the waning of his energy and the approach of death" (see \cite{Hes}, 143). Yet, one ceases to be enthusiastic about truth upon realizing how extraordinarily ill-equipped man is to recognize it, and that such a quest cannot constitute humanity's fundamental activity. Persisting in loyalty to this immature ideal leads to a weakening of our already limited capacity for spiritual awakening and the intuitive recognition of divine truth.

Admittedly, age brings sclerosis and reduced blood supply to the brain, resulting in a decline of our intellectual faculties. Yet "these afflictions ultimately have their upside," observes Hesse (see \cite{Hes}, 144), for we no longer perceive things with the same precision and intensity as before. In this sense, certain pains vanish, and one drifts toward a state of insensibility. Of all the accumulated knowledge and curiosities one might possess, only the ultimate fascination remains: that centered on the great, implacable figure of death. Consequently, the time for idle chatter or diversion has passed. "When one is old and has fulfilled one's task, one has the right to approach death in silence. Death has no need of mankind; it knows us -- it has seen enough of us. What it demands is silence. It is unseemly to disturb it, to speak to it, or to torment it with our idle talk. Upon reaching the threshold of its realm, we ought to pass by as though no one dwelt there" (see \cite{Hes}, 144). \section{Towards a Mathematical Formulation}

How can we mathematically formulate the new vision of old age that has just been proposed? A summary of the positive attributes of old age mentioned by previous authors might look like this:

A. A forward-looking and more or less prophetic sense; predominance of the spirit; quiet pleasures; proximity to death. 

B. Fusion with nature; spiritual and religious concerns; questioning the meaning of things and a quest for authenticity; mystical experiences of the world's mystery; patience and tranquility; access to beauty; a taste for repetition and routine; fascination with and a serious attitude toward death; the right to silence. 

In short:

1. Predominance of the forces of the spirit. 

2. Mystical experience, fusion with nature, and access to beauty. 

3. A taste for repetition and routine. 

4. Curiosity regarding death; 

5. Silence.

\subsection{Increased Spiritual Capacities and Gradient Dynamics}

It is evident that formalization based on declining functions can, at best, apply only to the physiology of aging. Moreover, the very notion of "decline" can be interpreted in a non-negative light -- much like Nietzsche’s reference to the "decline of Zarathustra," which mirrors the setting of the sun: for the mature individual, having passed the ideal -- and, so to speak, celestial -- phase of study, the task is to return to the earth and to real people. Old age could thus be interpreted in this sense: the decline in physiological functions might well correspond to an increased attentiveness to earthly matters and the concrete aspects of life. This suggests, rather, a gradient dynamic.

Recall that in mathematics and physics, the gradient of a function of several variables is a vector field that combines the various partial derivatives at each point, thereby indicating both the direction of the steepest local change and the magnitude of that change. The gradient of a function $f$ is denoted as $\text{grad}(f)$ or, using the "nabla" operator, as $\nabla f$ (sometimes with arrows placed above the symbol).

At every point where it is defined, the dot product with the gradient constitutes the function's differential -- that is, the linear part of its first-order Taylor expansion. This method allows a function of several variables to be locally approximated by a linear form, and the concept extends to real-valued functions defined on a Riemannian manifold.

It is worth noting that the gradient is always orthogonal to level lines or isosurfaces. It also enables the expression of constrained optimization conditions and plays a role in numerical analysis methods used to generate minimizing sequences.

Its use therefore seems entirely justified in the context of aging, where the goal is to optimize compromised living conditions while minimizing the most negative aspects of this process.

Let us immediately consider an example.

In a Euclidean Cartesian coordinate system, the gradient of a function $f$ that is differentiable at the point $a=(x_{1},x_{2},\dots ,x_{n})$ is the vector denoted by $\nabla f(a)$, with components $\tfrac {\partial f}{\partial x_{i}} (a)$ (where $i = 1, 2, ..., n$) -- that is, the partial derivatives of $f$ with respect to the coordinates at point $a$:
\begin{equation}
\nabla f(a)=
\begin{bmatrix}
\frac {\partial f}{\partial x_{1}}(a)\\\vdots \\ \frac {\partial f}{\partial x_{n}}(a)
\end{bmatrix}.
\end{equation}

\begin{figure}[h] %  figure placement: here, top, bottom, or page
\hspace{7\baselineskip}
\vspace{0\baselineskip}
\includegraphics[width=4in]{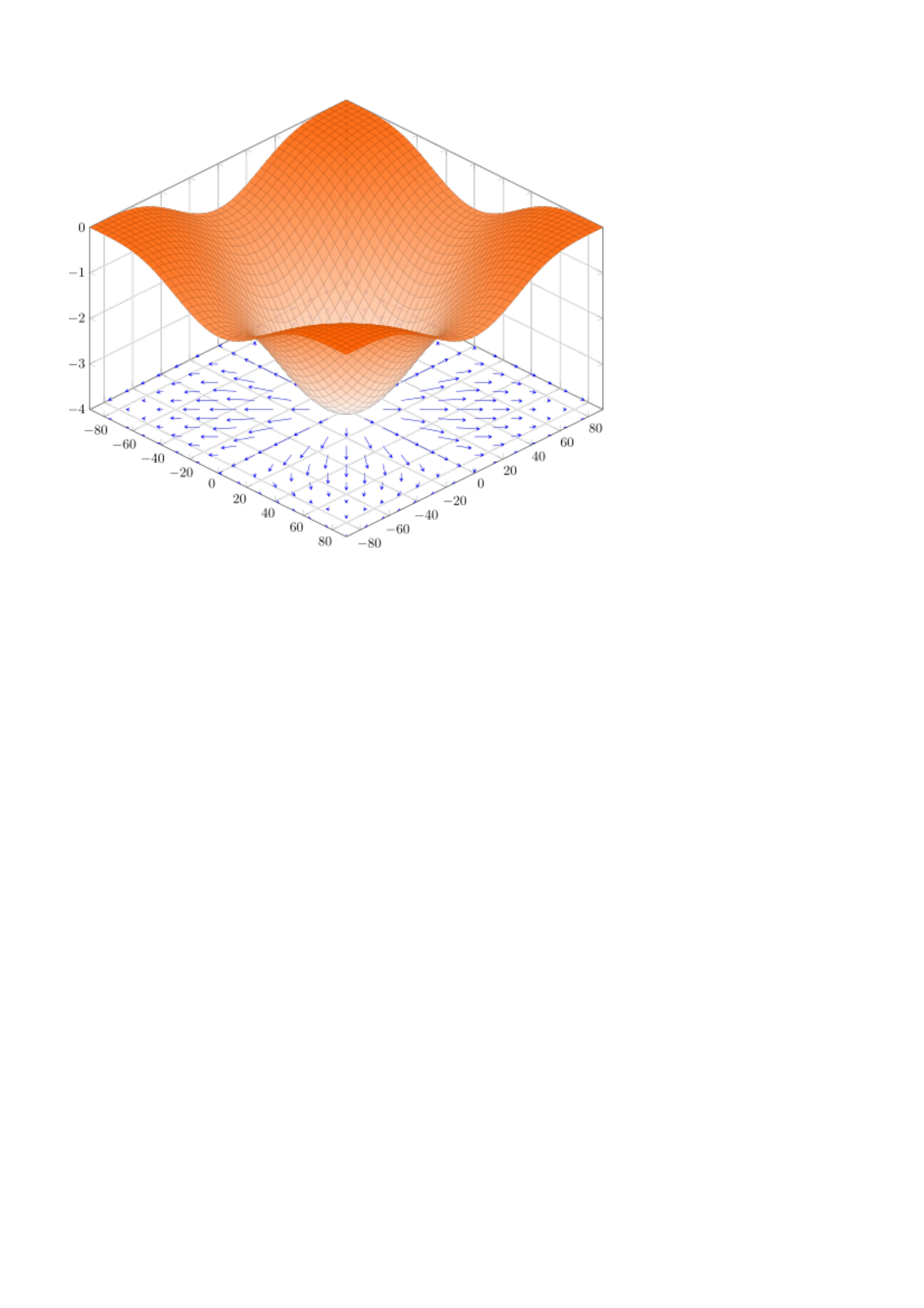}
\vspace{-14\baselineskip}
\caption{Scalar field (orange surface) and vector field (blue arrows)}
\label{fig: imm4}
\end{figure}

Thus, the scalar field of the function $f(x,y) = - (cos^2x + cos^2y)^2$ is represented by the orange surface in Fig. 4. In the same figure, the gradient of $f$ is a vector field represented by the blue arrows, each pointing in the direction where $f$ increases most rapidly.

In an orthonormal coordinate system, if the gradient vector is non-zero, it points in the direction of the function's steepest increase, and its magnitude equals the rate of increase in that direction.

The components of the gradient of $f$ are the coefficients of the variables in the equation of the tangent space to the graph of $f$ at point $a$. This property allows it to be defined independently of the choice of coordinate system, as a vector field whose components transform when switching from one coordinate system to another.

In this case, we have:
\[
\nabla f(x, y) = \left(\frac{\partial f}{\partial x}, \frac{\partial f}{\partial y}\right).
\]
Let:
\[
g(x) = \cos^2 x + \cos^2 y, \qquad \text{whence:} \quad f = -g^2,
\]
and a simple calculation yields:
\begin{equation}
\nabla f(x, y) = 2(\cos^2 x + \cos^2 y)(\sin(2x), \sin(2y)).
\end{equation}
Since $f \le 0$, its maxima are equal to 0 and are attained when $\cos x = \cos y = 0$, that is, when:
\[
x = \frac{\pi}{2} + k\pi, \qquad y = \frac{\pi}{2} + \ell\pi.
\]
At these points, the gradient is zero.

The minima are equal to -4 when $\cos^2 x = \cos^2 y = 1$, that is, for:
\[
x = k\pi \quad \text{and} \quad y = \ell\pi,
\]
and here, too, the gradient is zero.

If $f$ corresponds, for instance, to a potential of interest or attention directed toward certain concrete realities (people, things, knowledge, etc.), then $\nabla f$ indicates the direction of the steepest ascent of this potential, while $-\nabla f$ corresponds to the force of attraction or the intensity of the flow toward these realities. Formula (33) thus describes a field that is periodic in $x$ and $y$, featuring a structure of wells and peaks that repeat every $\pi$. This field can be viewed as a two-dimensional undulating landscape where trajectories follow paths of descent toward the minima $(k\pi, \ell\pi)$ before rising again -- a "bipolar" structure, in short. The generalization of the gradient to differentiable, vector-valued functions of several variables (as well as to differentiable mappings between Euclidean spaces) is the Jacobian matrix. The generalization to functions between Banach spaces is the Fréchet derivative.

\subsection{Mystical experience and mathematical structure}

What Hermann Hesse discusses oscillates between a dissolution of the boundary between the self and the world, a holistic view of reality -- or a unitary limit thereof -- and, finally, a vision of beauty that is simultaneously topological and poetic. Let us examine these various elements.

Let $S$ be the set representing the self and $N$ the set representing nature. Typically,
\begin{equation}
S \cap N \neq \emptyset,
\end{equation}
in other words, interactions exist, yet the two related entities remain distinct.

Conversely, mystical union could be represented by the equation $S \cup N = N$, or better yet, the equation:
\begin{equation}
S \subseteq N = N,
\end{equation}
which signifies that the self is no longer perceived as separate from the natural totality, or that it has completely dissolved into it.

It should be noted that equation (35) has exactly the same form as equation (11), of which it is, in a sense, a limit when one considers -- instead of two successive moments -- the totality of subjective as well as objective experience.

From a holistic perspective, one could also view each living being as a local function of a global system:
\begin{equation}
s_i = f_i(N),
\end{equation}
where $s_i$ represents an individual. In this case, the mystical experience would correspond to the realization that:
\begin{equation}
\forall s_i, s_i \subset N,
\end{equation}
with individual identity being merely a particular expression of the whole.

All of this could also be expressed as a limit approaching unity. Thus, if $d(S, N)$ measures the conceptual distance between self and nature, we would have:
\begin{equation}
\lim_{t \to \infty} d(S; N) = 0.
\end{equation}
In general terms, the mystical experience is linked to the disappearance of the illusion of separation between the self and the world -- a process that may be thought to unfold with the passing of time. In this sense, equation (35) remains the most elegant formulation of the phenomenon\footnote{This inclusion-equality can take on rather complex forms in a philosophy such as Spinoza's -- a philosophy whose simultaneously mathematical and mystical nature has been highlighted by one of its commentators, Hubbeling. For this, see our discussion and notes in \cite{Par2}, pp. 15–16.}.

\subsection{Repetitions, habits, and graphs with loops}

Let us now define a directed graph $G(V, E)$, where $V$ is the set of possible behavioral states and $E$ is the set of transitions between these states (waking up, breakfast, work, meals, etc.).

What is commonly referred to as a {\it habit} will be represented here as a cycle in the graph. Mathematically, a cycle is a sequence of vertices $(v_1, v_2, \dots, v_n, v_1)$ that returns to its starting point.

A habit thus corresponds to a trajectory that frequently traverses this cycle.

Since not all habits have the same strength, a weight $w_{ij} \ge 0$ can be assigned to each transition. The more the behavior is repeated, the more this weight increases. A deeply ingrained habit will therefore correspond to a heavily weighted edge. Aging is thus characterized not only by an increase in the cycles within an individual's behavioral graph but also by an increase in the weights associated with transitions between two nodes of the same cycle.

A more probabilistic representation of the situation would involve using the Markov chain formalism.

\begin{dfn}
A Markov chain (see Fig. 5) on $X$ with transition matrix $P$ is a sequence of random variables $(X_n)_{n \in \mathbb{N}}$ defined on a space ($\Omega, \beta, P$) and taking values ​​in $X$, such that, for all $n$ and all points $x_0,...,x_{n+1}$,
\[
\mathbb{P}[X_{n+1} = x_{n+1} | X_0 = x_0,...,X_n = x_n] = \mathbb{P}(x_n,x_{n+1}).
\]
Note that the conditional distribution $\mathbb{P}[X_{n+1} | (X_0,..., x_n)]$ is the transition probability $P(X_n,\cdot)$.
\end{dfn}

We have a transition matrix of the following form:
\[
\begin{pmatrix}
p_{11} & p_{12} &  \ldots \\
p_{21} & p_{22} & \ldots \\
\vdots & \vdots & \ddots
\end{pmatrix}
\]
where:
\[
p_{ij} =P(X_{t+1} = j | X_t = i)
\]
represents the probability of transitioning from one state to another. \begin{figure}[h] %  figure placement: here, top, bottom, or page
\hspace{7\baselineskip}
\vspace{0\baselineskip}
\includegraphics[width=4in]{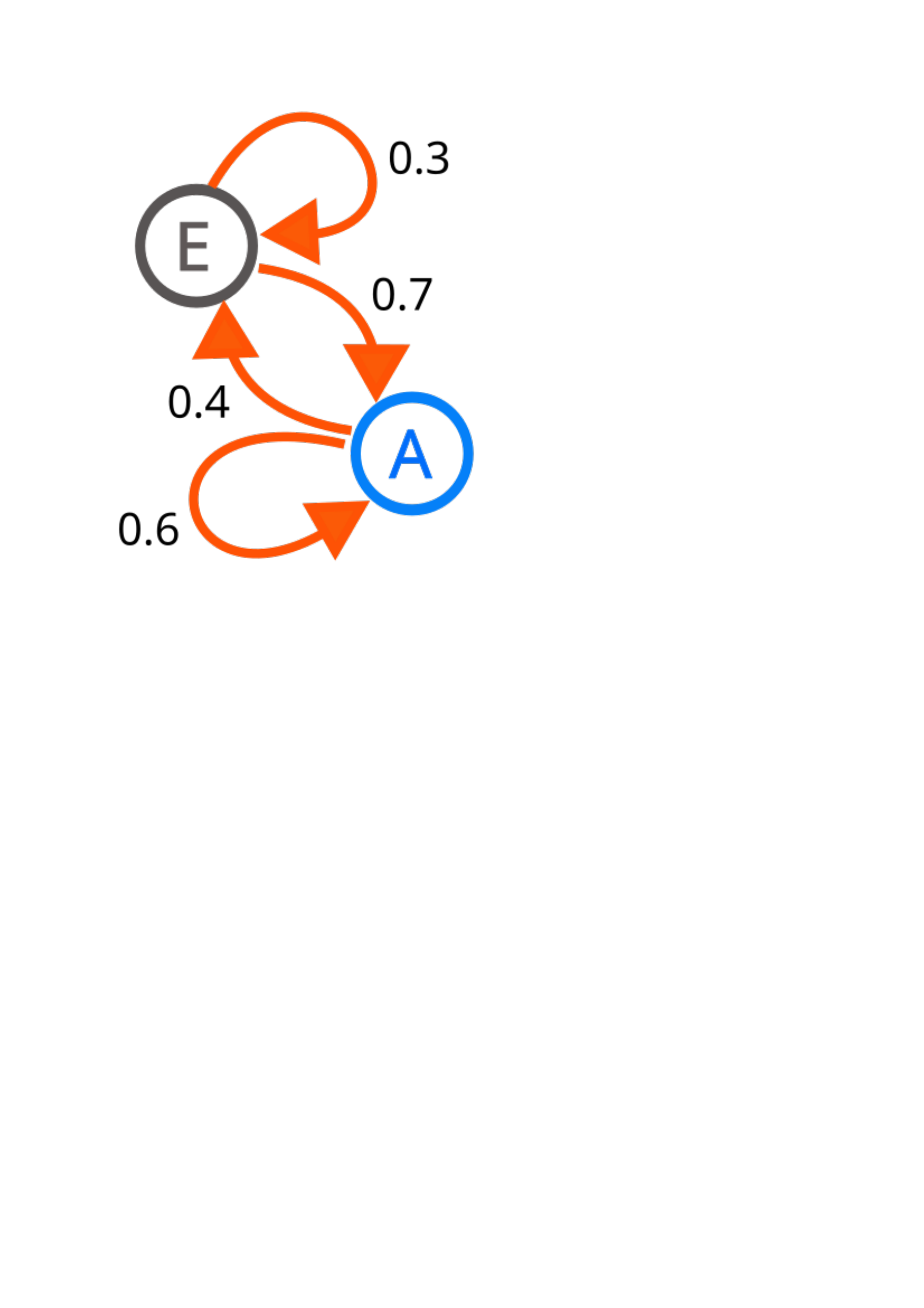}
\vspace{-14\baselineskip}
\caption{Example of a Markov chain}
\label{fig:  MM1}
\end{figure}
In this context, a habit corresponds to a region of the graph where return probabilities are high. The number of such regions would increase with age.

Another interesting formalization would be to view a habit as a dynamic attractor. In other words, given a sequence of behaviors $X_0, X_1, X_2, \dots$, certain configurations would particularly attract the trajectory. In the language of dynamical systems, the habit would be a behavioral {\it basin of attraction}.

One can also make the weights dependent on the number of repetitions:
\[
w_{ij}(t+1) = w_{ij}(t)  + \alpha N_{ij}(t),
\]
where:
$N_{ij}(t)$ is the number of times the transition has been performed, and $\alpha$ is a learning coefficient.

Thus, the more a loop is traversed, the easier it becomes to traverse, which aligns well with observations made in psychology and neuroscience.

It should be noted that complex habits can be represented as graphs of graphs, allowing for the creation of graph hierarchies such as $G_0 \subset G_1 \subset G_2 \dots$ etc., or even fractal structures.

\subsection{The pervasive presence of death in old age}

A human life corresponds to a finite number of events if one counts them in a rough manner (breaths, heartbeats, etc.). However, if one considers the potential density of events at any given moment (mathematical instants, thoughts, micro-events...), then the set of instants forms an uncountable continuum. It is infinite, albeit bounded in time. From birth to death -- excluding death itself -- it constitutes an interval $(0, T)$ homeomorphic to $\mathbb{R}$.

Bachelard wrote a remarkable text on this subject -- one generally misunderstood by philosophers, who lack sufficient mathematical background. Here it is in its entirety. "What gives duration its emotional character -- the joy or pain of being -- is the proportion or disproportion of the hours of life spent as hours of thought or hours of sympathy. Matter neglects to be; life neglects to live; the heart neglects to love. It is while sleeping that we lose Paradise. Let us follow the perspective of our own indolence: the atom radiates and frequently exists -- it utilizes a vast number of instants -- yet it does not utilize {\it all} instants. The living cell is already more sparing of its efforts; it utilizes only a fraction of the temporal possibilities afforded to it by the aggregate of atoms that constitute it. As for thought, it utilizes life only in irregular flashes. Three filters through which all too few instants reach consciousness! Thus, we feel a dull ache when we go in search of lost instants. We recall those rich hours marked by the thousand sounds of Easter bells -- those bells of resurrection whose tolls go uncounted because they all count, because each one finds an echo in our awakened soul. And this memory of joy turns to remorse when we compare those hours of total life with hours that are intellectually sluggish -- because they are relatively impoverished -- or with dead hours -- because they are empty (empty of purpose, as Carlyle said from the depths of his sorrow) -- or with hostile, interminable hours that yield nothing.
And we dream of a divine hour that would yield everything. Not merely the 'full' hour, but the 'complete' hour. The hour in which every instant of time would be utilized by matter; the hour in which every instant realized within matter would be utilized by life; the hour in which every living instant would be felt, loved, and thought. An hour, consequently, in which the relativity of consciousness would be erased, since consciousness would perfectly match the measure of time." ...complete.
Ultimately, objective time is maximal time; it is the time that encompasses all instants. It consists of the dense set of the Creator's acts" (see \cite{Bach} 47-48).

Obviously, human time is far removed from this dense set; it is, rather, a discrete subset -- though, if we adopt our initial postulate, a subset that remains infinite and can thus be equated with the countable, or in other words, with the set $\mathbb{N}$.

Consequently, one could view death as a compactification of life, treating it as a single point at infinity -- denoted by $\omega$ -- akin to the Alexandroff compactification.
This solution has the merit of aligning with one of the possible interpretations of the Schwartz space $\mathcal{S}$ of rapidly decreasing functions, which was introduced in our analytical formalization of aging. Indeed, this interpretation is noted by Schwartz himself in (\cite{Sch2}, 235). After demonstrating that the space $\mathcal{S}$ is locally convex, complete, and possesses a countable neighborhood basis -- meaning it is a Montel space where bounded sets and relatively compact sets coincide -- Schwartz proceeds to develop another possible representation of this space.

"Let us consider," he writes, "the $n$-dimensional sphere $S^n$, with the differentiable structure it possesses as a sphere within the $(n+1)$-dimensional Euclidean space $\mathbb{R}^{n+1}$. We know that $S^n$ can be represented by the space $\mathbb{R}^n$ augmented by a point at infinity, $\omega$; the differentiable structure is that of $\mathbb{R}^n$ at finite distances, while in the neighborhood of the point $\omega$, it is the structure obtained on $\mathbb{R}^n$ near the origin via inversion centered at the origin."
Since $S^n$ is compact, the space $\mathcal{D}_{S^n}$ of infinitely differentiable functions with compact support is simply the space of all infinitely differentiable functions on $S^n$. It possesses a countable neighborhood basis, and since $\mathbb{R}^n$ is an open subset of $S^n$, the restriction of any function $\bar{\phi}$ on $S^n$ is a function $\phi$ on $\mathbb{R}^n$. Through this correspondence between $\bar{\phi}$ and $\phi$, the topological space $\mathcal{S}$ on $\mathbb{R}^n$ is isomorphic to the closed vector subspace of $\mathcal{D}_{S^n}$ consisting of all functions $\bar{\phi}$ that -- along with all their successive derivatives -- vanish at the point $\omega$. This constitutes a topological isomorphism, allowing all properties of $\mathcal{S}$ to be derived from the properties of $\mathcal{D}_{S^n}$.

However, one might consider this solution too impoverished. Such a formalization would not allow for distinguishing between different ways of ending, nor for encoding the totality of life's asymptotic properties. For life is not merely a sequence -- convergent or otherwise. This is a set of events to which bounded continuous functions can be associated (for example: emotional intensity, moral value, impact on others...).

It would therefore be more fruitful to use the Stone-\v{C}ech compactification. Recall that if $X$ is a completely regular topological space, the Stone-\v{C}ech compactification is the largest compactification of $X$, in the sense that any continuous map from $X$ to a compact space extends uniquely to $\beta X$. Moreover, it is extremely rich, and its remainder $\beta X \smallsetminus X$ is enormous -- unless, of course, $X$ is already compact.

With the Stone-\v{C}ech compactification:

\begin{itemize}
\item Every bounded continuous function defined on life (on the interval $(0, T)$) extends continuously to $\beta(0, T)$.
\item Death, as a point (or set of points) in $\beta(0, T) \smallsetminus (0, T)$, encodes all possible limits of all these functions.
\item In other words, death contains the trace of all quantities that possess an asymptotic limit--even if, in real life, many of these limits do not exist. In fact, the Stone-Čech compactification forces them to exist.
\end{itemize}

On $\mathbb{N}$, the set of "human" moments, $\beta\mathbb{N}$ is interpreted as the set of ultrafilters. An ultrafilter on discretized life (days or key moments) represents a consistent way of deciding what holds true at the end. In a human life, different "ends" are possible, depending on the chosen criterion: an ultrafilter of happy moments, an ultrafilter of significant moments, or an ultrafilter of the memories one retains from one's past? In the Stone-\v{C}ech model, death is not a single point but an infinity of points (as many as there are ultrafilters). Each point corresponds to a possible way of "closing off life" -- that is, a different limit depending on the perspective (biological, psychological, social, spiritual, etc.) one takes on it.

The Stone-\v{C}ech model, which might be the one accessed by old age, captures the essential incompleteness of open life ($0, T$), without determining a single closure. Death, as $\beta(0, T) \smallsetminus(0, T)$, is a rich, uncountable object that preserves all asymptotic information. Of course, as a biological instant, it may appear to be just a point. But in reality, it is rather a horizon that contains all the possible endings of the story we have lived, and, in this sense, the Stone-\v{C}ech compact model models death as a concept, as a memorial and symbolic totality that keeps track of all ending scenarios. In this sense, it is larger than life, because $\beta X$ contains strictly more points than $X$. In this sense too, it surpasses life, completes it, transcends it by adding to it the infinity of possible limits. This is why only old age, which has lived and seen much, can approach its essence more closely than other periods of life.

\section{Conclusion}

We thus arrive at the idea that there is not just one mathematical model of old age. The analytic theory of declining functions provides only a superficial idea, which is certainly corrected by the algebraic approaches we have mentioned with idempotent (inclined) semirings or the age of relational structures. But with these structures, we remain primarily on the side of the negative. In reality, getting older brings compensatory elements: a gradient dynamic associated with spiritual functions previously neglected, possible access to a mystical experience, the development of new habits and pools of attraction reinforced by repetitions, and finally, access to an unprecedented vision of death, a multiple but necessary death, linked to the different ways of closing life.

\end{document}